\documentclass{article}
\usepackage[utf8]{inputenc}
\usepackage{graphics}
\usepackage{graphicx}
\usepackage{epsfig}
\usepackage{amssymb}
\usepackage[english,francais]{babel}
\newtheorem{theorem}{Theorem}[section]

\newtheorem{e-proposition}[theorem]{Proposition}

\newtheorem{e-definition}[theorem]{Definition\rm}

\newcommand{\R}{\mathbb{R}}

\renewcommand{\epsilon}{\varepsilon}
\newcommand{\eps}{\varepsilon}

\usepackage{amsmath}

\def\og{\leavevmode\raise.3ex\hbox{$\scriptscriptstyle\langle\!\langle$~}}
\def\fg{\leavevmode\raise.3ex\hbox{~$\!\scriptscriptstyle\,\rangle\!\rangle$}}

\title{A kinetic eikonal equation}

\author{Emeric Bouin\footnote{UMR CNRS 5669 'UMPA' and INRIA project 'NUMED', \'Ecole Normale Sup\'erieure de Lyon, 
46, all\'ee d'Italie, 
F-69364 Lyon Cedex 07, 
France. Mail : emeric.bouin@ens-lyon.fr, vincent.calvez@ens-lyon.fr }, Vincent Calvez $^*$ }

\begin{document}
\maketitle
\selectlanguage{english}
\begin{abstract}

%
We analyse the linear kinetic transport equation with a BGK relaxation operator. We study the large scale hyperbolic limit $(t,x)\to (t/\eps,x/\eps)$. We derive a new type of limiting Hamilton-Jacobi equation, which is  analogous to the classical eikonal equation derived from the heat equation with small diffusivity. We prove well-posedness of the phase problem and convergence towards the viscosity solution of the Hamilton-Jacobi equation. This is a preliminary work before analysing the propagation of reaction fronts in kinetic equations.


\vskip 0.5\baselineskip

\selectlanguage{francais}
\begin{center} \noindent{\bf R\'esum\'e} \vskip 0.5\baselineskip \noindent \end{center}
{\bf Une \'equation eikonale cin\'etique.}
Nous analysons une \'equation cin\'etique lin\'eaire de transport avec un op\'erateur de relaxation BGK. Nous \'etudions la limite hyperbolique de grande \'echelle $(t,x)\to (t/\eps,x/\eps)$. Nous obtenons \`a la limite une nouvelle \'equation de Hamilton-Jacobi, qui est l'analogue de l'\'equation eikonale classique obtenue \`a partir de l'\'equation de la chaleur avec petite diffusion. Nous d\'emontrons le caract\`ere bien pos\'e de l'\'equation v\'erifi\'ee par la phase, ainsi que la convergence vers une solution de viscosit\'e de l'\'equation de Hamilton-Jacobi. Ceci est un travail pr\'eliminaire en vue d'analyser la propagation de fronts de r\'eaction pour des \'equations cin\'etiques.
  


\end{abstract}

\selectlanguage{francais}
\section*{Version fran\c{c}aise abr\'eg\'ee}

Nous  consid\'erons  un mod\`ele cin\'etique lin\'eaire avec un op\'erateur de relaxation BGK, pos\'e sur un ensemble de vitesses $V$ sym\'etrique et born\'e. On analyse le comportement de l'\'equation dans la limite hyperbolique de grande \'echelle $(t,x)\to \left(\frac t \eps,\frac x \eps\right)$, 
\begin{equation} \label{KinSharp-fr}
\partial_t f^{\eps} + v\cdot \nabla_{x} f^{\eps} = \frac{1}{\epsilon} \left( M(v) \rho^{\eps}- f^{\eps} \right ), \quad (t,x,v) \in \R_+ \times \R^n \times V\, .
\end{equation}
Nous d\'emontrons que la phase $\varphi^\eps$ d\'efinie par la relation $f^{\eps}(t,x,v) = M(v) e^{- \frac{\varphi^{\eps}(t,x,v)}{\eps}}$ converge (localement) uniform\'ement, lorsque $\eps \to 0$, vers une fonction $\varphi^0(t,x)$ ind\'ependante de $v$. De surcro\^it, la fonction $\varphi^0$ est solution de viscosit\'e de l'\'equation de Hamilton-Jacobi suivante,
\begin{equation}\label{disp-fr}
\int_{V} \frac{M(v)}{1 - \partial_t \varphi^0(t,x) - v\cdot \nabla_{x} \varphi^0(t,x) } dv =1 \,, \quad (t,x) \in \R_+ \times \R^n \,.
\end{equation}
Cette \'equation peut se r\'e\'ecrire sous la forme canonique $\partial_t \varphi^0 + H(\nabla_x \varphi^0) = 0$ pour un hamiltonien effectif $H(p)$ qui est lipschitzien et convexe.  En r\`egle g\'en\'erale, nos travaux consistent \`a {\em homog\'en\'eiser} l'\'equation \eqref{KinSharp-fr} par rapport \`a la variable de vitesse. Le probl\`eme aux valeurs propres dans la {\em cellule} $V$ s'\'ecrit comme suit~: trouver un vecteur propre $Q(v)$ et une valeur propre $H(p)$ tels que
\begin{equation*}
\left( 1 + H(p) - v\cdot p  \right) Q(v) = \int_{V} M(v') Q(v') dv'\, .
\end{equation*}

La d\'emonstration du passage \`a la limite de \eqref{KinSharp-fr} vers \eqref{disp-fr} s'appuie sur une s\'erie d'estimations {\em a priori} qui d\'emontre que $\varphi^\eps$ appartient \`a l'espace de Sobolev $W^{1,\infty}$, avec un contr\^ole uniforme en $\eps>0$ (Proposition \ref{solphi} ci-dessous). Dans un deuxi\`eme temps, nous d\'emontrons que toute fonction test $\psi^0(t,x)$ de classe $\mathcal C^2$ telle que $\varphi^0 - \psi^0$ admet un maximum local en  $\left( t^0, x^0 \right)$, v\'erifie 
\begin{equation*}
\int_{V} \frac{M(v)}{1 - \partial_t \psi^0(t^0, x^0) - v \cdot \nabla_{x} \psi^0(t^0, x^0) } dv \leq 1\,.
\end{equation*} 
Ceci d\'emontre que $\varphi^0$ est une sous-solution de viscosit\'e l'\'equation de Hamilton-Jacobi \eqref{disp-fr}. Un raisonnement identique montre qu'il s'agit aussi d'une sur-solution de viscosit\'e. La d\'emonstration se base sur la construction d'un correcteur microscopique $\eta(t,x,v)$ d\'efini de fa\c{c}on ad-hoc par la relation
\begin{equation*}
\forall (v,v')\in V\times V\, , \quad e^{\eta(t,x,v)} - e^{\eta(t,x,v')} = \left( v' - v \right) \cdot  \nabla_{x}  \psi^0(t,x) \, .
\end{equation*}


\selectlanguage{english}
\section{Large-scale limit and derivation of the Hamilton-Jacobi equation}

We consider the following kinetic equation with BGK relaxation operator: 
\begin{equation} \label{Kin}
\partial_t f + v\cdot \nabla_{x} f = M(v) \rho - f\,, \quad (t,x,v) \in \R_+ \times \R^n \times V\,,
\end{equation} 
where $f(t,x,v)$ denotes the density of particles moving with speed $v \in V$ at time $t$ and position $x$. The function $\rho(t,x)$ denotes the macroscopic density of particules: 
\begin{equation*}
\rho(t,x) = \int_{V} f(t,x,v) \,dv\,, \quad (t,x) \in \R_+ \times \R^n\,.
\end{equation*}
Here $V$ denotes a bounded symmetric subset of $\R^n$. We assume that the Maxwellian $M$ is symmetric and satisfies the following moment identities: 
\begin{equation*}
\int_{V} M(v)dv = 1\,, \qquad \int_{V} v M(v) dv = 0\,, \qquad \int_{V}v^2 M(v) dv = \theta^2\,.
\end{equation*}
In this paper we focus on the large scale  hyperbolic  limit  $(t,x) \rightarrow \left( \frac{t}{\eps} , \frac{x}{\eps} \right)$, $\eps \to 0$. The kinetic equation \eqref{Kin} reads as follows in the new scaling:
\begin{equation} \label{KinSharp}
\partial_t f^{\eps} + v\cdot \nabla_{x} f^{\eps} = \frac{1}{\epsilon} \left( M(v) \rho^{\eps}- f^{\eps} \right ), \quad (t,x,v) \in \R_+ \times \R^n \times V\, .
\end{equation}
Clearly, the velocity distribution relaxes rapidly towards the Maxwellian distribution. This motivates the introduction of the following Hopf-Cole transformation: 
\begin{equation*}
\label{HC}
f^{\eps}(t,x,v) = M(v) e^{- \frac{\varphi^{\eps}(t,x,v)}{\eps}}\,.
\end{equation*}
where we expect the phase $\varphi^{\eps}$ to become independent of $v$ as $\eps \to 0$. To avoid technical complications due to ill-prepared data, we set $\varphi^{\eps}(0,x,v) = \varphi_{0}(x)\geq 0$ as an initial data for \eqref{KinSharp}.
The equation satisfied by $\varphi^{\eps}$ reads:
\begin{equation}\label{eqphi}
\partial_t \varphi^{\eps} + v\cdot \nabla_{x} \varphi^{\eps} = \int_V M(v') \left( 1 -  e^{\frac{\varphi^{\eps} - \varphi^{\eps\prime}}{\eps}} \right) dv'\,, \quad (t,x,v) \in \R_+ \times \R^n \times V\,,
\end{equation}


\begin{theorem}\label{VS}
Let $V \subset \R^n$ be bounded and symmetric, and $M \in L^1(V)$ be nonnegative and symmetric. Then $\varphi^{\eps}$ converges (locally) uniformly towards $\varphi^0$, where $\varphi^0$ does not depend on $v$. Moreover $\varphi^0$ is the viscosity solution of the following Hamilton-Jacobi equation:
\begin{equation}\label{disp}
\int_{V} \frac{M(v)}{1 - \partial_t \varphi^0(t,x) - v\cdot \nabla_{x} \varphi^0(t,x) } dv =1 \,, \quad (t,x) \in \R_+ \times \R^n \,.
\end{equation}
The denominator of the integrand is positive for all $v\in V$.
\end{theorem}

The last observation in Theorem \ref{VS} is not compatible with an unbounded velocity set. 

One can understand this as an homogenization problem in the velocity variable. Moreover, the associated eigenproblem in the cell $V$ writes: Find an eigenvector $Q(v)$ such that: 
\begin{equation*}
\label{eigenpb}
\left( 1 - \partial_t \varphi^0 - v\cdot \nabla_{x} \varphi^0  \right) Q(v) = \int_{V} M(v') Q(v') dv'\, .
\end{equation*} 
This eigenproblem can be solved explicitly, and yields formula \eqref{disp}.

Thanks to monotonicity properties, 
we can boil down to the classical framework of first order Hamilton-Jacobi equations. Indeed, writing equation \eqref{disp} as $G(\partial_t \varphi^0, \nabla_x \varphi^0) = 0$, we observe that $G$ is increasing with respect to the first variable. Hence the equation is equivalent to $\partial_t \varphi^0 + H(\nabla_x \varphi^0) = 0$, where the effective Hamiltonian $H$ is defined through the implicit formula, 
\begin{equation}\label{eq:effective H}
\int_V \frac{M(v)}{(1 + H(p) - v\cdot p)}\, dv = 1\, .
\end{equation}
Differentiating \eqref{eq:effective H} we obtain,
\[ \int_V \frac{M(v)}{(1 + H(p) - v\cdot p)^2} \left( \nabla H(p) - v \right)\, dv = 0\, . \]
We deduce $\|\nabla H\|_\infty \leq V_{\rm max}$. This is in accordance with the underlying kinetic equation, since $ \nabla  H$ can be interpreted as the group speed, which is bounded by the maximal speed of the particles. 
Differentiating \eqref{eq:effective H} twice we obtain
\[ \left(\int_V \frac{M(v)}{(1 + H(p) - v\cdot p)^2}   \, dv\right)  \mathrm{D}^2 H(p) = 2 \int_V \frac{M(v)}{(1 + H(p) - v\cdot p)^3} \left( \nabla H(p) - v \right)\otimes \left( \nabla H(p) - v \right)\, dv  \, . \]
We deduce that the effective Hamiltonian is convex.

As an example, we can compute the effective Hamiltonian $ H $ in one dimension for a constant Maxwellian $M \equiv \frac12$ on $V = ( -1,1 )$. We obtain $ H(p) = \frac{p - \tanh(p)}{\tanh(p)}$. It is equivalent to $\theta^2 |p|^2$ for small $p$ ($\theta^2 = \frac13$). Another example where the effective hamitonian is explicit is given by the Maxwellian $M(v) = \frac{1}{2} \left( \delta_1 + \delta_{-1} \right)$, though it is not a $L^1$ function. This corresponds to a two velocities model (also known as the {\em telegraph equation}, see \cite{Fedotov-1999,Bouin-Calvez-Nadin-2011}). In this case we obtain the relativistic hamiltonian $H(p) = \frac{\sqrt{1 + 4p^2} - 1 }{2}$.

Interestingly enough, we obtain a Hamilton-Jacobi equation which differs from the classical eikonal equation. The latter could have been expected from the following argumentation. 
The formal limit of equation \eqref{KinSharp} at order $O(\eps)$ is the heat equation with small diffusivity:
\begin{equation*}
\partial_t \rho^\eps = \eps \theta^2 \Delta_x \rho^\eps\,, \quad (t,x) \in \R_+ \times \R^n\, .
\end{equation*}
It is well-known that the phase $\phi^\eps = -\eps \log \rho^\eps$ satisfies in the limit  $\eps \to 0$ the classical eikonal equation in the sense of viscosity solutions \cite{Evans-Ishii-1985,Freidlin-1986,Evans-1989,Evans-Souganidis-1989,Freidlin-Wentzell}:  
\begin{equation}\label{eik}
\partial_t \phi^0 + \theta^2 \vert \nabla_x \phi^0 \vert^2 = 0\, .
\end{equation}
We only have asymptotic equivalence between the two approaches for small $|p|$ as can be seen directly on  \eqref{eq:effective H} by Taylor expansion: $H(p)\sim \theta^2|p|^2$. 

%

\begin{figure}
\begin{center}
\includegraphics[width = .48\linewidth]{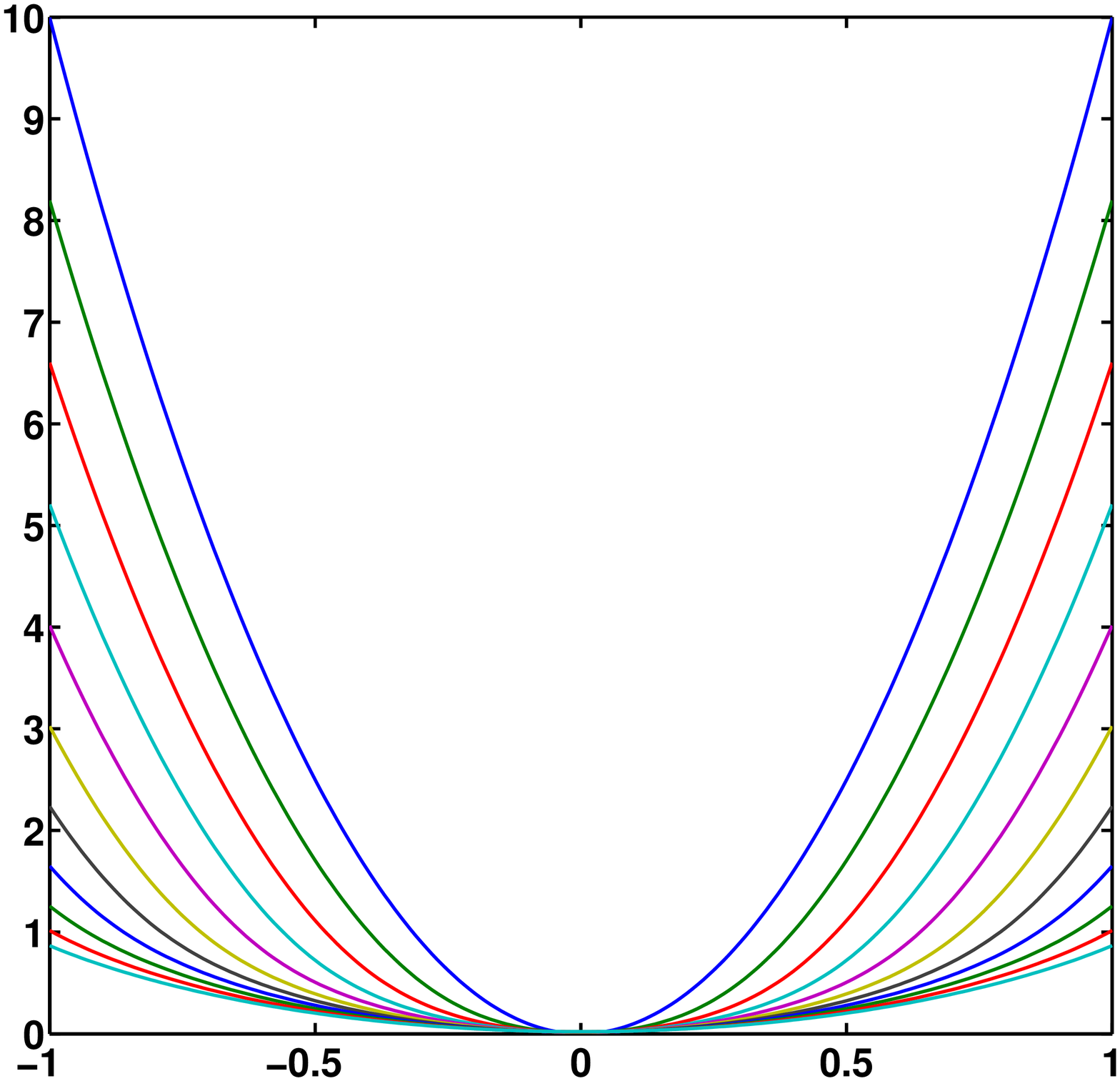} \,
\includegraphics[width = .48\linewidth]{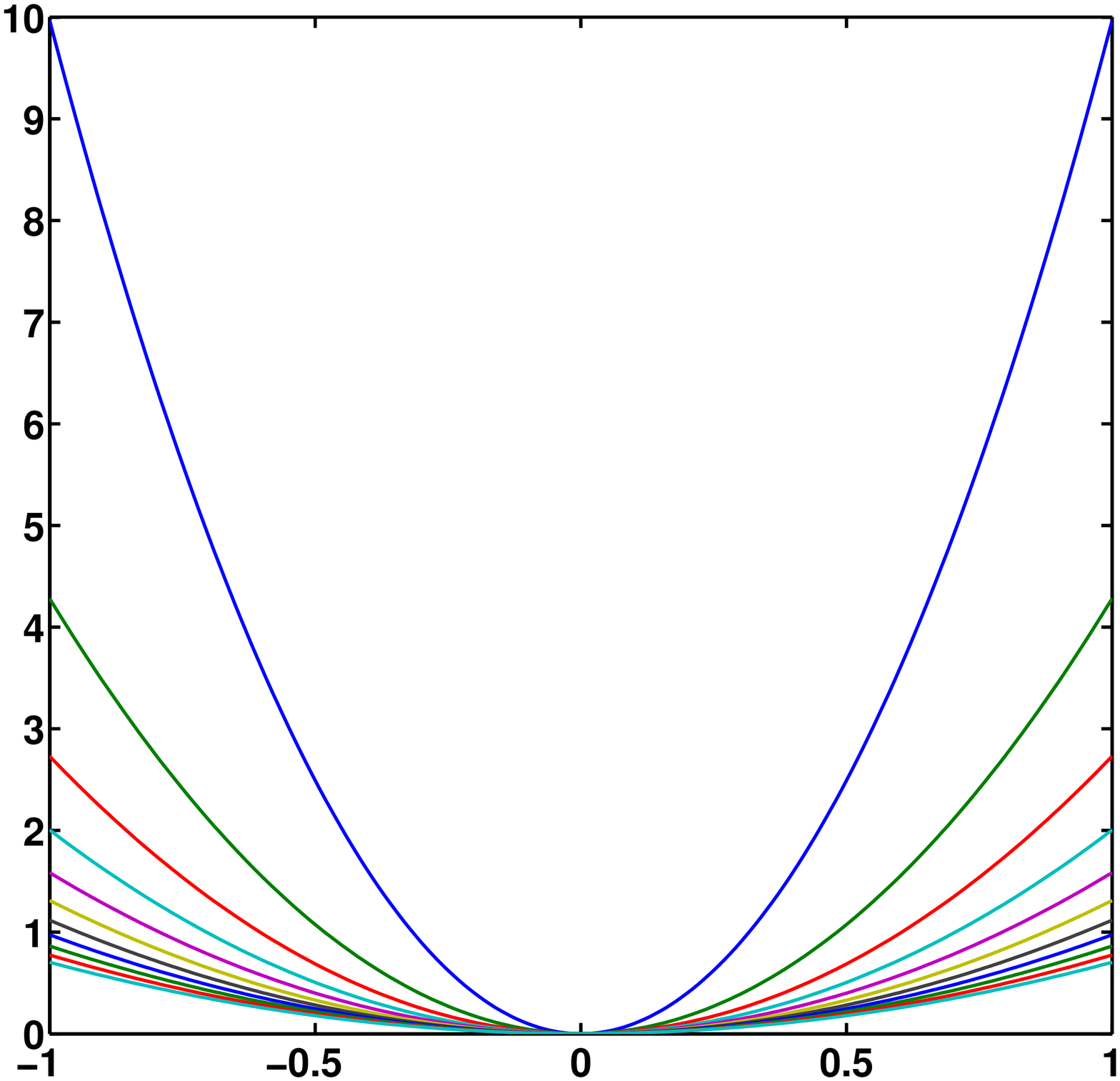} \\ 
\caption{Numerical simulations of the Hamilton-Jacobi equation $\partial_t \varphi + H(\partial_x \varphi) = 0$. ({\em left}) The kinetic eikonal equation \eqref{disp} where $M(v) = \frac12 \mathbf{1}_{(-1,1)}$. ({\em right}) The classical eikonal equation $H(p) = \theta^2 |p|^2$ \eqref{eik}.  In both cases the initial condition is a parabola.}
\label{fig:comparison}
\end{center}
\end{figure}

In Figure \ref{fig:comparison} we show numerical simulations of the kinetic eikonal equation \eqref{disp}, with a constant Maxwellian on $V = ( -1,1 )$, and we compare it with the classical eikonal equation \eqref{eik}.  


We end this introduction by listing some possible extensions of Theorem \ref{VS} for other choices of transport and scattering operators. We will develop a more general framework in a future work.
\begin{enumerate}
\item In the case of the Vlasov-Fokker-Planck equation, 
\begin{equation*}
\partial_t f^\eps + v \cdot \nabla_{x} f^\eps - \nabla_x V(x) \cdot \nabla_v f^\eps =  \frac1\eps \nabla_v \cdot \left( \nabla_v f^\eps + v f^\eps \right)\,, \quad (t,x,v) \in \R_+ \times \R^n \times \R^n\,, 
\end{equation*}
we obtain simply the eikonal equation $\partial_t \phi^0 + \vert \nabla_x \phi^0 \vert^2 = 0$ in the WKB expansion $f^\eps = M(v)e^{-\frac{\varphi^\eps}\eps}$, where $M(v)$ is a Gaussian. 
\item It is challenging to replace the BGK operator in \eqref{KinSharp} by a convolution operator $L(f) = K * f - f$, where $K$ is a probability kernel \cite{Barles-Mirrahimi-Perthame-2009}. However in this case we are not able to solve explicitly the eigenproblem in the cell $V$. 
\item In a forthcoming work we will investigate the propagation of reaction fronts in kinetic equations, following \cite{Evans-Souganidis-1989,Fedotov-1999,Cuesta-Hittmeir-Schmeiser-2010,Bouin-Calvez-Nadin-2011}.
\end{enumerate}

\section{Proof of Theorem \ref{VS}}

First let us mention that the solution $\varphi^\eps$ remains nonnegative for all times. We proceed in two steps. First we prove uniform estimates with respect to $\eps>0$. It allows to extract a uniformly converging subsequence. Second we identify the limit as the viscosity solution of equation \eqref{disp} using the maximum principle. The second step relies on the construction of a suitable corrector $\eta(t,x,v)$ \cite{Evans-1989,Evans-2004}.
%
\medskip

\noindent\textbf{Step 1. Existence and uniform bounds.}\medskip

\begin{e-proposition}\label{solphi}
Let $V\subset \R^n$ be a bounded subset. Assume $M \in L^1(V)$ and $\varphi_0 \in W^{1,\infty} \left( \R^n \right)$. The kinetic equation \eqref{eqphi} has a unique solution $\varphi^{\eps} \in W^{1,\infty}\left( \R_+ \times \R^n \times V \right)$. Furthermore, the solution satisfies the following uniform estimates:
\begin{align}
\label{bound}
& 0 \leq \varphi^{\eps}(t, \cdot ) \leq \Vert \varphi_0 \Vert_{\infty}\,,\\
\label{dxbound}
& \Vert \nabla_{x} \varphi^{\eps}(t, \cdot ) \Vert_{\infty} \leq \Vert \nabla_{x} \varphi_0 \Vert_{\infty}\,, \\
\label{dvbound}
 &\Vert \nabla_v \varphi^{\eps}(t, \cdot ) \Vert_{\infty} \leq t  \Vert \nabla_x \varphi_0 \Vert_{\infty}   \, , \\
\label{dtbound}
& \Vert \partial_t \varphi^{\eps}(t, \cdot ) \Vert_{\infty} \leq  V_{\rm max}\Vert  \nabla_x \varphi_0 \Vert_{\infty}\, .
\end{align}
\end{e-proposition}

\noindent \textbf{Proof.}
We obtain a unique solution $\varphi^\eps$ from a fixed point method on the Duhamel formulation of \eqref{eqphi}:
\begin{equation}\label{eq:duhamel}
\varphi^{\eps}(t,x,v) = \varphi_0(x - tv ) + \int_{0}^{t} \int_{V} M(v') \left( 1 - e^{\frac{  \varphi^{\eps}(t-s,x - sv,v) - \varphi^{\eps}(t-s,x-sv,v') }{\eps}} \right)  dv' ds,
\end{equation} 
We obtain directly, 
\begin{equation*}
\forall \eps > 0,  \qquad 0 \leq \varphi^{\eps}(t,x,v) \leq \varphi_0(x - tv ) + t\, .
\end{equation*} 
This ensures that $\varphi^{\eps}$ is uniformly bounded on $[0,T] \times \R^n \times V$.
To prove the bound \eqref{bound}, we define $\psi_{\delta}^{\eps} (t,x,v) = \varphi^{\eps}(t,x,v) - \delta t - \delta^4 \vert x \vert^2$. For any $\delta > 0$, $\psi_{\delta}^{\eps}$ attains a maximum at point $(t_{\delta},x_{\delta},v_{\delta})$. Suppose that $ t_{\delta} > 0$. Then, we have 
\begin{equation*} 
\partial_t \varphi^{\eps} (t_{\delta},x_{\delta},v_{\delta}) \geq \delta, \qquad \nabla_{x} \varphi^{\eps} (t_{\delta},x_{\delta},v_{\delta}) = 2\delta^4 x_{\delta}.
\end{equation*}
As a consequence, we have at the maximum point $(t_{\delta},x_{\delta},v_{\delta})$:
\begin{equation}\label{ppmax}
0\geq \int_V M(v') \left( 1 -  e^{\frac{    \psi_{\delta}^{\eps} (t_{\delta},x_{\delta},v_{\delta}) - \psi_{\delta}^{\eps} (t_{\delta},x_{\delta},v')  }{\eps}} \right) dv' \geq \delta + 2 v_{\delta} \delta^4 x_{\delta} \geq \delta - 2V_{\rm max}\delta^4 \vert x_{\delta} \vert\,.
\end{equation}
Moreover, the maximal property of $(t_{\delta},x_{\delta},v_{\delta})$  also implies 
\begin{equation*}
\Vert  \varphi^{\eps} \Vert_{\infty} - \delta^4 \vert x_{\delta} \vert^2  \geq \varphi^{\eps}(t_{\delta},x_{\delta},v_{\delta}) - \delta t_{\delta} - \delta^4 \vert x_{\delta} \vert^2 \geq \varphi^{\eps}(0,0,v_{\delta}) \geq 0 \, .
\end{equation*} 
We obtain a contradiction with \eqref{ppmax} since $\delta^{-6}/(2V_{\rm max}) \leq \vert x_{\delta} \vert^2 \leq   \delta^{-4} \Vert  \varphi^{\eps} \Vert_{\infty}$ cannot hold for sufficiently small $\delta>0$. As a consequence  $t_{\delta} = 0$, and we have,
\begin{equation*}
\forall (t,x,v) \in [0,T] \times \R^n \times V, \qquad \varphi^{\eps}(t,x,v) \leq \varphi^0(x_{\delta},v_{\delta}) + \delta t + \delta^4 \vert x \vert^2 \leq \Vert \varphi_0 \Vert_{\infty} + \delta t + \delta^4 \vert x \vert^2.
\end{equation*}
Passing to the limit $\delta \to 0$, we obtain \eqref{bound}.
To find the bound \eqref{dxbound}, we  use the same ideas on the difference $\varphi_h^{\eps} (t,x,v) = \varphi^{\eps}(t,x+h,v) - \varphi^{\eps}(t,x,v)$. The equation for $\varphi_h^{\eps}$ reads as follows,
\begin{equation*}
\partial_t \varphi_h^{\eps}  + v \cdot \nabla_{x} \varphi_h^{\eps} = \int_V M(v') e^{\frac{\varphi^{\eps} - \varphi^{\eps\prime}}{\eps}} \left( 1 -  e^{\frac{    \varphi_h^{\eps} - \varphi_h^{\eps\prime}  }{\eps}} \right) dv'.
\end{equation*}
Using the same argument as above with a $\delta-$correction, we conclude that
\begin{equation*}
\forall (t,x,v) \in [0,T] \times \R^n \times V, \qquad \varphi_h^{\eps}(t,x,v) \leq \sup_{(x,v) \in \R \times V}\left \vert  \varphi^0(x+h,v) - \varphi^0(x,v) \right \vert
\end{equation*}
The same argument applies to $- \varphi_h^{\eps}$,
\begin{equation*}
\partial_t \left( - \varphi_h^{\eps} \right) + v \cdot \nabla_{x} \left( - \varphi_h^{\eps} \right) = - \int_V M(v') e^{\frac{\varphi^{\eps} - \varphi^{\eps\prime}}{\eps}} \left( 1 -  e^{- \frac{    \left( - \varphi_h^{\eps} \right) - \left( - \varphi_h^{\eps\prime} \right) }{\eps}} \right) dv'.
\end{equation*}
so that the r.h.s has the right sign when $- \varphi_h^{\eps}$ attains a maximum.
Finally, 
\begin{equation*}
\forall (t,x,v) \in [0,T] \times \R^n \times V, \qquad \vert \varphi_h^{\eps}(t,x,v) \vert \leq \sup_{(x,v) \in \R \times V}\left \vert  \varphi^0(x+h,v) - \varphi^0(x,v) \right \vert \leq \left \Vert \nabla_x \varphi^0 \right\Vert_{\infty} \vert h \vert.
\end{equation*}
from which the estimate \eqref{dxbound} follows.

To obtain regularity in the velocity variable \eqref{dvbound}, we differentiate \eqref{eqphi} with respect to $v$, 
\begin{equation*}
\left( \partial_t + v \cdot \nabla_x \right) \left( \nabla_v \varphi^{\eps} \right) = - g_{\eps}(\varphi^{\eps}) \nabla_v \varphi^{\eps} - \nabla_x \varphi^{\eps},
\end{equation*}
where $g_{\eps}(\varphi^{\eps}) = \frac{1}{\eps} \int_{V} M(v') e^{ \frac{ \varphi^{\eps} - \varphi^{\eps\prime}   }{\eps} } dv' \geq 0 $. Multiplying by $\frac{\nabla_v \varphi^\eps}{|\nabla_v \varphi^\eps|}$, we obtain 
\begin{align}
\left( \partial_t + v \cdot \nabla_x \right) \left( \vert \nabla_v \varphi^{\eps} \vert  \right) & = -  g_{\eps}(\varphi^{\eps}) \vert \nabla_v \varphi^{\eps} \vert  -   \left( \nabla_x \varphi^{\eps} \cdot \frac{\nabla_v \varphi^\eps}{|\nabla_v \varphi^\eps|} \right) \label{eq:nabla v phi} \\
& \leq   \Vert \nabla_x \varphi_{0} \Vert_{\infty} \nonumber \,.
\end{align}
from which we deduce \eqref{dvbound} since  $\nabla_v\varphi_0 = 0$ by hypothesis. 

Finally, the bound \eqref{dtbound} is obtained similarly as the bound on $\nabla_{x} \varphi^{\eps}$ \eqref{dxbound}, using the difference $\varphi_s^{\eps} (t,x,v) = \varphi^{\eps}(t+s,x,v) - \varphi^{\eps}(t,x,v)$. We obtain 
\begin{equation*}
\forall (t,x,v) \in [0,T] \times \R^n \times V, \qquad \vert \varphi_s^{\eps}(t,x,v) \vert \leq \sup_{(x,v) \in \R \times V}\left \vert  \varphi^\eps(s,x,v) - \varphi^0(x,v) \right \vert \, .
\end{equation*}
We use the Duhamel formulation \eqref{eq:duhamel} to estimate the last contribution:
\begin{equation*}
\left \vert  \varphi^\eps(s,x,v) - \varphi^0(x,v) \right \vert \leq \left\vert  \varphi_0(x-sv ) - \varphi_0(x ) \right \vert + o(s)\, .
\end{equation*}
The estimate \eqref{dtbound} follows.
$\hfill \square$

%
%
%

\medskip

\noindent\textbf{Step 2. Viscosity solution procedure.}\medskip

\noindent From Proposition \ref{solphi} we deduce that the familly $\left( \varphi^{\eps} \right)_{\eps}$ is locally uniformly bounded in $W^{1,\infty}\left( \R_+ \times \R^n \times V \right)$. 
Then, from the Ascoli-Arzel\`a theorem, we can extract a locally uniformly converging subsequence. We denote by $\varphi^0$ the limit. Furthermore, from the fact that $\int_V M(v') e^{\frac{\varphi^{\eps} - \varphi^{\eps\prime}}{\eps}} dv'$ is uniformly bounded on $[0,T]\times \R^n\times V$, we deduce that $\varphi^0$ does not depend on $v$.

Let $\psi^0 \in \mathcal{C}^{2}\left( \R_+ \times \R^n \right)$ be a test function such that $\varphi^0 - \psi^0$ has a local maximum at  $\left( t^0, x^0 \right)$. We want to show that $\psi^0$ is a subsolution of \eqref{disp}, yielding that $\varphi^0$ is a  viscosity subsolution \cite{Crandall-Evans-Lions-1984}. The supersolution case can be performed similarly. Thereby, we define a corrective term $\eta$ not depending on $\eps$:  $\psi^{\eps} = \psi^0 + \eps \eta$. 
The corrector $\eta$ is defined up to an additive constant. We choose the renormalization $\int_{V} M(v') e^{- \eta'} dv' = 1$. We define $\eta$ as follows,
\begin{equation}\label{eta}
\forall (v,v')\in V\times V\, , \quad e^{\eta(t,x,v)} - e^{\eta(t,x,v')} = \left( v' - v \right) \cdot  \nabla_{x}  \psi^0(t,x) \, .
\end{equation} 
The corrector $\eta$ is well defined. In fact, we can choose any $v_0\in V$ and define $e^{\eta(t,x,v)}  = \mu_0 +  \left( v_0 - v \right) \cdot  \nabla_{x}  \psi^0(t,x)$. There is a unique positive $\mu_0 = e^{\eta(t,x,v_0)}$ under the condition $\int_{V} M(v') e^{- \eta'} dv' = 1$.

The  uniform convergence ensures that $\varphi^{\eps} - \psi^{\eps}$ has a maximum at $(t^\eps,x^{\eps},v^{\eps})$, where $(t^\eps,x^{\eps})$ is close to $(t^0,x^0)$. As $V$ is a bounded set, the sequence $( v^{\eps} )$ has an accumulation point, say $v^*$. 
We can extract a subsequence (without relabelling) such that $(t^{\eps },x^{\eps },v^{\eps }) \to (t^0,x^{0}, v^{*})$. 
We have at $(t^{\eps},x^{\eps},v^{\eps})$:
\begin{equation*}
1 - \partial_t \psi^{\eps} - v^{\eps} \cdot \nabla_{x} \psi^{\eps} = 1 - \partial_t \varphi^{\eps} - v^{\eps} \cdot \nabla_{x} \varphi^{\eps} = \displaystyle \int_{V} M(v') e^{\frac{\varphi^{\eps} - \varphi^{\eps\prime}}{\eps}} dv' \, .
\end{equation*}   
From the maximum property of $(t^\eps,x^{\eps},v^{\eps})$, the last inequality implies at this point :  
\begin{equation*}
1 - \partial_t \psi^{\eps }(t^\eps,x^{\eps},v^{\eps}) - v^{\eps } \cdot \nabla_{x} \psi^{\eps }(t^\eps,x^{\eps},v^{\eps}) \geq \int_{V} M(v') e^{\eta(t^\eps,x^{\eps},v^{\eps})  - \eta(t^\eps,x^{\eps},v^{\prime})} dv'.
\end{equation*}
Passing to the limit, we obtain at $(t^0, x^0)$:
\begin{equation*}
1 - \partial_t \psi^0(t^0, x^0) -  v^*   \cdot \nabla_{x} \psi^0(t^0, x^0) \geq  \int_{V} M(v') e^{\eta(t^0, x^0,v^*)- \eta(t^0, x^0,v')} dv' = e^{\eta(t^0, x^0,v^*)} \,.
\end{equation*}  
From the very definition of the corrector $\eta$ \eqref{eta}, this writes also:
\begin{equation*}
\forall v\in V\,, \quad 1 - \partial_t \psi^0(t^0, x^0) -  v   \cdot \nabla_{x} \psi^0(t^0, x^0) \geq    e^{\eta(t^0, x^0,v)} \,.
\end{equation*}  
Therefore we obtain at point $(t^0, x^0)$,
\begin{equation*}
\int_{V} \frac{M(v)}{1 - \partial_t \psi^0(t^0, x^0) - v \cdot \nabla_{x} \psi^0(t^0, x^0) } dv \leq \int_{V} M(v ) e^{- \eta(t^0, x^0,v)} dv = 1\,.
\end{equation*} 
We conclude that $\psi^0$ is a subsolution of  \eqref{disp}. 

%
$\hfill \square$






\end{document}